\documentclass[12pt]{amsart}

\usepackage{amsmath,amsthm,amscd,euscript}
\def \r{\mathbb R}
\def \q{\mathbb Q}
\def \z{\mathbb Z}

\newtheorem{theorem}{Theorem}[section]

\newtheorem{statement}{Statement}[section]

\newtheorem{hyp}{Conjecture}
\newtheorem{prob}{Problem}
\theoremstyle{remark}

\theoremstyle{plane}

\theoremstyle{definition}
\newtheorem{definition}{Definition}[section]
\newtheorem{example}{Example}[section]

\title
[Two-dimensional continued fractions.]
{On examples of two-dimensional periodic continued fractions.}
\author{Oleg Karpenkov}
\date{3 November 2004}
\email[Oleg Karpenkov]{karpenk@mccme.ru}

\begin{document}
\input epsf
\maketitle
\tableofcontents

\section*{Introduction.}

This paper is a short survey of the recent results on examples of
periodic two-dimensional continued fractions (in  Klein's model).
These examples were constructed in the last fifteen years by many
mathematicians in their works. In the last part of this paper we
formulate some questions, problems and conjectures on geometrical
properties concerning to this subject.

Klein's model of the multidimensional continued fractions is one
of the most well-known generalizations of one-dimensional
continued fractions to the multidimensional case. The definition
of Klein continued multidimensional fractions was given by
F.~Klein in his works~\cite{Kl1} and~\cite{Kl2} (later on we will
call them just multidimensional continued fractions). F.~Klein
considered the integer approximations for zeroes of quadratic
forms $f(x,y)=ax^2+bxy+cy^2$, where $f(x,y)$ is a product of two
linear factors. Ordinary continued fractions provide a good
algorithm for these approximations. Suppose the set of all zeroes
for the form $f(x,y)$ does not contain any integer point other
than the origin, and the linear factors of the form $f(x,y)$ are
not proportional. Consider the convex hulls of all integer points
in each of the four connected components of the complement to the
set of zeroes of the form $f(x,y)$. The boundaries of the convex
hulls are called sails and the union of them --- the
one-dimensional continued fraction (the notions of ordinary
continued fraction and one-dimensional continued fraction are
differs a little, the relation between these notions is described
in Section~1). It turns out that the vertices of the sails provide
the best in some sense integer approximations of zeroes for the
form $f$ (the detailed description see in~\cite{Hin}). The main
idea of the algorithm is to construct the vertices of the sail
recursively. The description of the algorithm is not in the
limelight of the present work, so we omit it (it can be found,
for example, in~\cite{Arn2}).

The idea of F.~Klein was to find some effective generalization of the algorithm
to the case of homogeneous cubic forms of three variables
that is decomposable in a product of three linear factors.
He generalized the one-dimensional continued fractions to the
two-dimensional case.
F.~Klein considered the convex hulls of all integer points in each of the eight
connected components of the complement to the set of all zeroes for the cubic
form (we again suppose that three planes of zeroes are in general position
and there is no integer points in this planes other
than the origin).
The boundaries of the convex hulls are called the sails and the union
of the convex hulls --- the two-dimensional continued fraction.
This definition generalizes to the case the other dimensions
(see Section~1 for the precise definitions)
Unfortunately, F.~Klein failed to find an effective algorithm that
generalizes the one-dimensional algorithm.

The interest for multidimensional continued fractions was revived by
V.~I.~Arnold's work~\cite{Arn1}.
He realized the reasons of Klein's failure and proposed a new approach
to the study of the multidimensional continued fractions.
We describe one of in a few words.
The one-dimensional sails are infinite polygonal lines and hence
the set of all vertices possesses the natural ordering.
This ordering is used in the algorithm.
The multidimensional sails are infinite polyhedrons without
any natural ordering of the vertices.
In order to overcome such difficulties
V.~I.Arnold suggested to investigate the geometry and combinatorics
of the sails at first (it contains the study
of properties for such notions as
affine types of the faces of the sails,
their quantities and frequencies,
integer angles between the faces, integer distances, volumes
and so on, see also~\cite{Arn2}).

Klein multidimensional continued fractions have many relations
with some other bran\-ches of mathematics. For example,
J.-O.~Moussafir showed the connection between the sails of multidimensional
continued fractions and Hilbert bases
in his work~\cite{Mus}.
In~\cite{Tsu} H.~Tsuchihashi found the relation between periodic
multidimensional continued fractions and multidimensional cusp singularities
which generalizes the relation between
the ordinary continued fractions and two-dimensional cusp singularities.

We will mostly consider periodic Klein continued fractions
(i.e. the fractions that are invariant under some particular
group action described in Section~1) in this paper.
Some important properties concerning to non-periodic Klein
multidimensional continued fractions are described in the
works of G.~Lachaud~\cite{Lac4}, J.-O.~Moussafir~\cite{Mus} and
O.~N.~German~\cite{Ger}.
The other generalizations of continued fractions can
be found in the work~\cite{Min} of H.~Minkowski and the work~\cite{Voro}
of G.~F.~Voronoi. Some attempts of
comparison of the efficiency for the
algorithms for calculation different types
of continued fractions are given in~\cite{BP2} and~\cite{BP3}
by A.~D.~Bruno and V.~I.~Parusnikov.

This paper is organized as follows.
In the first section we give all definitions and notions that are
necessary to our subject.
In the second section we show the simplest three examples of
torus decompositions for cubic irrationalities.
Then we present the classification of all continued
fractions of cubic irrationalities associated to the operator with the
small norm (the sum of absolute numbers of the coefficients is less
or equal than six).
In the third section we give some examples of the infinite sequences
of nonequivalent torus decompositions.
Finally, in the fourth section we formulate some interesting problems
and pose appropriate conjectures.

{\bf Acknowledgement.} The author is grateful to professor V.~I.~Arnold for constant attention
to this work and for good remarks and questions, to N.~G.~Moshevitin
and R.~Uribe-Vargas for useful discussions
and CEREMADE (Paris-Dauphine) for the hospitality.

\section{Notions and definitions.}

In this section we recall some basic notions and definitions
(see also~\cite{Kar1}).
Consider the space $\r^{n+1}$ ($n\ge 1$)
over the field of real numbers.
The point of $\r^{n+1}$ is said to be {\it integer} if
all coordinates of this point are integer.
The polyhedron is said to be {\it integer} if
all its vertices are integer.
Consider a set of $n+1$ hyperplanes of $\r^{n+1}$ passing through the origin
in general position.
The complement to the union of these hyperplanes consists of
$2^{n+1}$ open orthants. Let us choose an arbitrary orthant.
\begin{definition}
The boundary of the convex hull of all integer points
except the origin
in the closure of the orthant is called {\it the sail}.
The union of all $2^{n+1}$ sails of the space $\r^{n+1}$ is called
{\it $n$-dimensional continued fraction} associated to
the given $n+1$ hyperplanes in general position in
$(n+1)$-dimensional space.
\end{definition}

Two $n$-dimensional sails are said to be {\it equivalent}
if there exists a linear transformation that preserves the integer lattice
of the $(n+1)$-dimensional space and maps one continued fraction to the other.

Two $n$-dimensional continued fractions are said to be {\it equivalent}
if there exists a linear transformation that preserves the integer lattice
of the $(n+1)$-dimensional space and maps one continued fraction to the other.

The definition of one-dimensional continued fraction is a little
different to the definition of ordinary continued fraction.
We will briefly discuss the relation between these two definitions.
The definition of the ordinary continued fraction for some positive real
number $a$ is given by the following construction.
Consider a ray $\{(x,y)|y=ax,x>0\}$ on the plane $OXY$
with the vertex in the origin. This ray divides the first quadrant
of the plane $Q=\{ (x,y)|x>1,y>1 \}$ into two parts. Let us construct two
boundaries of the convex hulls for the integer points in the closure
of each part without the origin.
These boundaries are polygonal lines. The polygonal line for each
boundary consists of a ray and an infinite numbers of segments
in the case of irrational number $a$,
and it consists of two rays and some finite number of segments
in the case of rational number $a$.
Consider the union of all segments of these two polygonal lines.
This set is called {\it the ordinary
continued fraction} for the positive real number $a$.

Consider some arbitrary one-dimensional continued fraction $L$
constructed by two lines passing through the origin.
Let $y=a_1x$ and $y=a_2x$ be the equations for these lines
(note that we regards now only positive irrational real
numbers $a_1$ and $a_2$).
The following statement holds.

\begin{statement}
There exist such radius $R>0$ that the union of ordinary continued
fractions for the numbers $a_1$ and $a_2$
coincides $($as a set$)$ with the one-dimensional continued fraction $L$
in the intersection of the quadrant $Q$ and the complement
to the disk with center in the origin and
radius equals $R$.
\end{statement}

The case of negative $a_1$ or $a_2$ are similar
and are left to the reader as an exercise.

Now let $A \in GL(n+1,\z)$ be the operator with irreducible
characteristic polynomial over the field $\q$ whose roots are all
real and distinct. Let us take the $n$-dimensional spaces that
spans all possible subsets of $n$ linearly independent
eigenvectors of the operator $A$. The spans of every $n$
eigenvectors defines $n+1$ hyperplanes passing through the origin
in general position. These hyperplanes defines the
multidimensional continued fraction.

\begin{definition}
The $n$-dimensional continued fraction associated to
the operator $A \in GL(n+1,\z)$ as above is called
{\it the $n$-dimensional continued fraction of $(n+1)$-algebraic
irrationality}.
The case of $n=1 (2)$ corresponds to {\it one$($two$)$-dimensional
continued fractions of quadratic $($cubic$)$ irrationalities}.
\end{definition}

The following statement holds (for proof see~\cite{Kar1}).

\begin{statement}
Continued fractions associated to the operators $A$ and $B$ of the group
$GL(n+1,\z)$ with distinct real irrational eigenvalues
are equivalent iff there exists such an integer operator $X$
with determinant equals one that the operator $\tilde A$, obtained
from the operator $A$ by means of the conjugation by the operator $X$
commutes with $B$.
\end{statement}

Now let us formulate the notion of periodic continued
fractions associated to some algebraic irrationalities.
Suppose the characteristic polynomial of the operator $A$
with all real roots be irreducible over the field of rational numbers.
Consider the group of integer operators with the determinant equals one
and positive eigenvalues
that commutes with $A$.
These operators form an Abelian group.
It follows from Dirichlet unity elements theorem (see.~\cite{BSh})
that this group is isomorphic to $\z^n$ and its action is free.
Each operator of this group preserves the integer lattice and
 the union of all $n+1$ hyperplanes,
and hence it preserves the $n$-dimensional continued fraction
(as far as all eigenvalues are positive the sails maps to themselves
also).
The factor of a sail under such group action is isomorphic to
$n$-dimensional torus at that.
The converse is also true:

\begin{statement}
For any double-periodic sail of a continued fraction
there exist an operator $A$ of the group $GL(n+1,\z)$
with irreducible over the field of rational numbers characteristic polynomial
with three real roots
such that one of the sails of the continued fraction associated to $A$
coincides with the given one.
\end{statement}

These two statements generalize Lagrange theorem for the ordinary
continued fractions
(see~\cite{Kor1},~\cite{Tsu},~\cite{Lac1} and~\cite{Lac3}).

Suppose some characteristic (object) is invariant under the action
of the group $SL(n+1,\z)$. Then this characteristic (object) is invariant
under the action of the subgroup that consists of all
integer operators with the determinant
equals one and positive eigenvalues that commutes with $A$.
Hence this characteristic (object) naturally corresponds to the induced
characteristic (object) for the factorization
of the sail by the described above group action.

The first example of such object is $k$-dimensional planes in $\r^{n+1}$
for some fixed $k\le n$.
Under the action of some operator from the group $SL(n+1,\z)$
any $k$-dimensional plane maps to some $k$-dimensional plane.
The second important example is sublattices of all integer points
in the $k$-dimensional subspaces.
Under the action of some operator from the group $SL(n+1,\z)$
any $k$-dimensional plane maps to some $k$-dimensional plane
and the integer sublattice of the first one maps to the integer sublattice
of the second plane in the one-to-one way.
To be more precise, the $k$-dimensional planes with the integer sublattice
are invariant under the group of all affine transformations
that preserves the integer lattice of the space $\r^{n+1}$.
This group contains the group $SL(n+1,\z)$ as a subgroup.

\begin{definition}
We say that two faces in the space $\r^{k}$ equipped with lattice
have the same {\it integer affine type} if
there exist a lattice preserving affine transformation of the space
that maps one of the faces to the other.
\end{definition}

It follows from this definition that integer affine type of the face is
invariant under the action of the group of all affine transformations,
that preserves the lattice of integer points.
As we have already mentioned
this group contains the group $SL(n+1,\z)$ as a subgroup.
Hence our tori are equipped with face decomposition.
All faces are some polyhedrons of the dimension less than or equal $n$.

Let us describe one important series of $SL(n+1,\z)$-action invariants.

\begin{definition}
The ratio of the Euclidean distance from an integer $k$-dimensional plane
(containing a $k$-dimensional integer lattice) to the integer point
to the minimal Euclidean distance from this plane to some
integer point in the complement of the plane
is called the {\it integer distance} from this plane to the
integer point.
\end{definition}

Let us rewrite the integer distance from the integer
$k$-dimensional plane to the integer point in the integer affine
invariant terms. Consider the lattice of integer points in the
$(k+1)$-dimensional plane that spans the given integer
$k$-dimensional plane and given integer point. Consider the
sublattice of this lattice generated by the given integer point
and all integer points in the given $k$-dimensional plane. Let us
choose some zero and standard addition for the lattice and the
sublattice. We have an Abelian group $\z^{k+1}$ and some its
subgroup. The index of this subgroup is equal to the integer
distance from the integer $k$-dimensional plane to the integer
point. This index is an integer affine invariant. So this
distance is also an integer affine invariant and therefore it is
a characteristic of the torus face decomposition.

Let us generalize the notion of distance to the case
of two integer planes of dimension $k_1$ and $k_2$
(we only need that each of them is not contained in the other one).
Consider the lattice of the integer points in the plane that spans
the given two planes.
Now consider the sublattice of this lattice generated by the
all integer points in the given two planes.

\begin{definition}
If the intersection of the given two planes is empty,
then the index of this sublattice is called {\it the integer distance}
between the planes.
If the intersection of the given two planes is not empty,
then corresponding index is called {\it the integer angle}
between the planes.
\end{definition}

All these indices are integer affine invariant,
and therefore all integer distances and integer angles are
also characteristic of the torus face decomposition.

{\it Remark.}
The last definition has natural generalization to the case
of $l$ different planes of dimensions $k_1, \ldots , k_l$.

\begin{definition}
We call by {\it torus decomposition corresponding to
the sail of the continued $n$-dimensional continued fraction
for the $(n+1)$-algebraic irrationality}
(or torus decomposition, for short) the above torus face decomposition
equipped with integer affine structure for all faces of all dimensions,
and all integer distances and integer angles
between the planes containing two faces, for any two faces.
\end{definition}

By {\it a fundamental region} of the sail we call a union of some
faces that contains exactly one face from each equivalence class.

\section{Three single examples of torus decompositions.}

Let us start this section with three single examples of the torus
decomposition for two-dimensional continued fractions of cubic
irrationalities. One can find a nice and rich collection of
beautiful single decompositions, and also a lot of useful
references concerning to the examples in the work~\cite{site} by
K.~Briggs. I hope that new examples from the present paper would
sooner or later appear in the list of the work of K.~Briggs.
 Further we work in $\r^3$. In this and the next sections we
describe only one of the eight possible sails for each continued
fraction.

To describe more precisely the torus decompositions for two-dimensional
continued fractions we need the following
geometric $SL(n+1,\z)$-invariants:
an integer length of a segment, an integer area of the face,
an integer angle between two rays, and an integer distance
from the plane (of a face) to some integer point.
The definition of the integer distance was given in the last part of the
previous section.
Let us give all other definitions.
They are all included in the definition of an
affine type of the face and hence are
$SL(3,\z)$-invariant.

\begin{definition}
The number of inner integer points of an integer segment in $\r^{3}$
plus one is called the {\it integer length} of the segment.
Integer points divide this segment in the congruent parts.
The ratio of the Euclidean area for an integer polygon in $\r^3$
to the Euclidean area for the minimal integer triangle lying in the
plane of the polygon is called its {\it integer area}.
\end{definition}

\begin{definition}
By {\it the integer angle} between two integer rays (i.e. containing
more than one integer point) with the vertex at the same integer point
we call the value
$S(u,v)/(|u|\cdot|v|)$, where $u$ and $v$ are arbitrary
integer vectors passing along the rays
and $S(u,v)$ is the integer volume of the triangle with edges $u$ and $v$.
\end{definition}

{\it Remark.}
The smallest possible integer volume of a triangle is equal to one
(in standard parallelogram measuring all area will be two times less).
The integer area of the triangle is equal to the
index of the lattice subgroup generated by its edges having the common vertex.
The integer angle between two integer rays with vertex in the point
is equal to the index of the lattice subgroup generated by all integer
points lying on edges of the angle.

\begin{example}\label{ex1}
Now all is done to describe the simplest example. We consider
two-dimensional fraction that generalizes the ordinary continued
fraction corresponding to the golden ratio
$\frac{1+\sqrt{5}}{2}$ $($see~\cite{Kor3}$)$.
All four one-dimensional sails of the one-dimensional continued
fraction corresponding to the golden ratio are equivalent.
They consist of the segments
of integer length one. The corresponding circle decompositions
consist of one vertex and one edge.
This fraction is given by the following operator:
$$
\left(
\begin{array}{cc}
1&1\\
1&2\\
\end{array}
\right)
$$

The generalization of this one-dimensional
continued fraction to the multidimensional case was
given by E.~I.~Korkina in the work~\cite{Kor2}.
The multidimensional continued fraction associated to the operator
$$
\left(
\begin{array}{ccccc}
1&1&1&\cdots&1\\
1&2&2&\cdots&2\\
1&2&3&\cdots&3\\
\vdots&\vdots&\vdots&\ddots&\vdots\\
1&2&3&\cdots&n\\
\end{array}
\right)
$$
is called {\it the generalized golden ratio}.

Here we consider the continued fraction associated to the generalized
gold ratio operator:
$$
\left(
\begin{array}{ccc}
1&1&1\\
1&2&2\\
1&2&3\\
\end{array}
\right).
$$

The torus decomposition corresponding to this operator is homeomorphic
to the following one:
$$
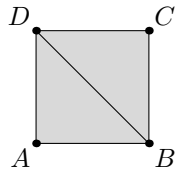
\epsfbox{fig3.1}$$
Here the segment $AB$ is identified with the segment
$DC$ and the segment $AD$ to the segment $BC$.

{\it Remark.} Note that only homeomorphic type of the faces is given
on the picture above.
This means, that we erase all integer points.
We have only the torus, equipped with the graph.
Each face of the graph corresponds to some face of the torus decomposition
of the considering periodic sail,
each segment of the graph corresponds to some segment of the torus decomposition,
and each vertex of the graph corresponds to some vertex of the torus
decomposition.
The integer affine types of the corresponding faces are given below.
Let us find all integer invariants mentioned above.

The integer areas of this two triangles are equal to~1.

The integer lengths of all three segments are equal to~1.

All six integer angles are equal to~1.

All two triangles have integer affine types of the simplest triangle:
$$
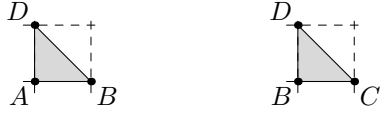
\epsfbox{fig3.2}$$

The integer distance from the triangle $ABD$ to the origin is equal to 2,
and from the triangle $BCD$ is equal  to~1.

This torus decomposition was constructed by
E.~I.~Korkina~\cite{Kor2}, G.~Lachaud~\cite{Lac2}
and A.~D.~Bryuno and V.~I.~Parusnikov~\cite{BP}
approximately at the same time.
\end{example}

\begin{example}\label{ex2}
The second example was given by A.~D.~Bryuno and
V.~I.~Parusnikov~\cite{BP}.
They construct the continued fraction that associated to
the following operator:
$$
\left(
\begin{array}{ccc}
1&1&1\\
1&-1&0\\
1&0&0\\
\end{array}
\right).
$$

The torus decomposition corresponding to this operator is also homeomorphic
to the following one:
$$
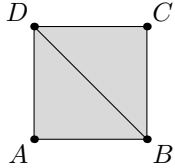
\epsfbox{fig3.1}$$
Here the segment $AB$ is identified with the segment
$DC$ and the segment $AD$ to the segment $BC$.

The integer area of the triangle $ABD$ is equal to~1 and
the integer area of the triangle $BCD$ is equal to 3.

The integer lengths of all three segments are equal to~1.

All three integer angles of the triangle $ABD$ are equal to~1.
All three integer angles of the triangle $BCD$ are equal to 3.

The integer distance from the triangle $ABD$ to the origin is equal to 2,
and from the triangle $BCD$ is equal  to~1.

The triangle $ABD$ has an integer affine type of the simplest triangle,
The triangle $BCD$ has an integer affine type of the triangle with the
vertices $(-1,-1)$, $(0,1)$, and $(1,0)$.
$$
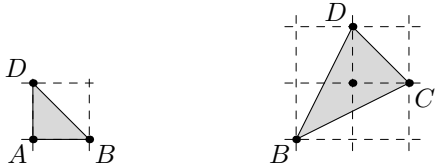
\epsfbox{fig3.3}$$
\end{example}

\begin{example}\label{ex3}
The third example was given by V.~I.~Parusnikov~\cite{Par1}.
This continued fraction is associated to
the following operator:
$$
\left(
\begin{array}{ccc}
0&1&0\\
0&0&1\\
1&1&-3\\
\end{array}
\right).
$$

The torus decomposition corresponding to this operator is also homeomorphic
to the following one:
$$
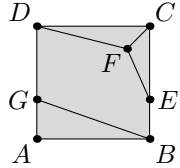
\epsfbox{fig3.4}$$
Here the segment $AB$ is identified with the segment
$DC$ and the polygonal line $AGD$ --- to the polygonal line $BEC$
(the point $G$ is identified with the point $E$).

The integer areas of all triangles are equal $1$ and
the integer area of the pentagon $BEFDG$ is equal to $5$.

The integer lengths of all three segments are equal to~1.

All integer angles of all triangle are equal to~1.
The angle in the vertex $D$ is equal to~1 and in other four vertices
the angles are equal to 2.

The integer distance from the triangle $ABG$  and the pentagon
$BEFDG$ to the origin is equal to~1,
from the triangle $CDF$ to the origin is equal to 2,
and from the triangle $CFE$ is equal  to 3.

All triangles have an integer affine type of the simplest triangle,
The pentagon $BEFDG$ has an integer affine type of the pentagon with the
vertices $(-1,0)$, $(-1,1)$, $(0,1)$ $(1,0)$ and $(1,-1)$:
$$
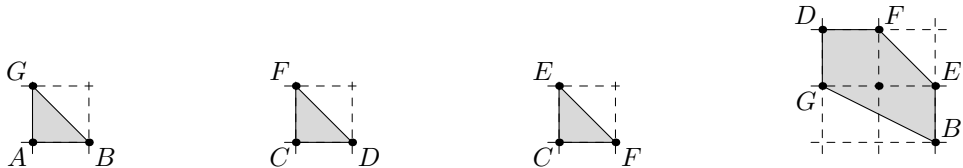
\epsfbox{fig3.5}$$
\end{example}

Some other single torus decomposition for some two-dimensional
continued fractions was investigated in the works~\cite{Par2} and~\cite{Par3}.

The continued fractions constructed in the examples~\ref{ex1},
\ref{ex2} and~\ref{ex3} are also known as the continued fractions
corresponding to the first, the second and the third Davenport form.

The first steps in classification of continued fractions with a small norm
based on these three examples were made by the author
of this issue in the work~\cite{Kar2}.
Let us now consider the following norm for the space of matrices:
the sum of absolute values of all the coefficients for the matrix.
\begin{theorem}
There is no hyperbolic matrix
among the set of matrices with integer coefficients and irrational eigenvalues,
and hence there is no continued fraction associated to the operator
with such norm.

If the norm of the hyperbolic matrix with
integer coefficients and irrational eigenvectors is equal to five
(their number is 48),
then the corresponding continued fraction is equivalent to the
generalization of the ordinary fraction for the golden ratio.
This fraction is shown in Example~\ref{ex1}.

The amount of such matrices with norm equals six is 912,
and only three different two-dimensional continued fractions associated
to them:
480 continued fractions are equivalent to the fraction
of Example~\ref{ex1},
240 continued fractions are equivalent to the fraction
of Example~\ref{ex3}
and 192 --- to the fraction of Example~\ref{ex2}
(for more
information see~\cite{Kar2}).
\end{theorem}

The classification of two-dimensional continued fractions
with the norm equals seven or more is unknown.

\section{Infinite sequences of torus decomposition examples.}

Now we continue the description of examples with the infinite
sequences of the torus decompositions.
The first two infinite sequences of torus decompositions was calculated by
E.~I.~Korkina in the work~\cite{Kor2}. One of this sequences is shown below.

\begin{example}\label{ex_k}
This sequence of continued fractions is associated to
the following operators for $m\ge 0$:
$$
\left(
\begin{array}{ccc}
0&0&1\\
1&0&-m-5\\
0&1&m+6\\
\end{array}
\right).
$$

The torus decomposition corresponding to this operator is homeomorphic
to the following one:
$$
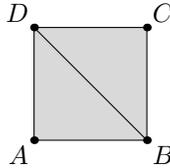
\epsfbox{fig3.1}$$
Here the segment $AB$ is identified with the segment
$DC$ and the segment $AD$ to the segment $BC$.

The integer areas of this two triangles are equal to~1.

The integer lengths of all three segments are equal to~1.

All six integer angles are equal to~1.

All two triangles have integer affine types of the simplest triangle:
$$
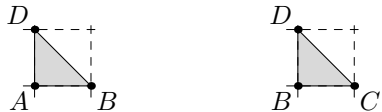
\epsfbox{fig3.2}$$

The integer distance from the triangle $ABD$ to the origin is equal to~$m+2$,
and from the triangle $BCD$ is equal  to~1.
\end{example}

Many other examples of infinite sequence of continued fractions
were done by the author of this issue in the work~\cite{Kar1}.
Now we show some of them.
The following sequence generalize the one from the previous example.

\begin{example}
This sequence depends on two integer parameters $a,b\ge 0$.
The continued
fractions of the sequence are associated to the following operators:
$$
\left(
\begin{array}{ccc}
0&1&0\\
0&0&1\\
1&1+a-b&-(a+2)(b+1)\\
\end{array}
\right).
$$

The torus decomposition corresponding to this operator is homeomorphic
to the following one:
$$
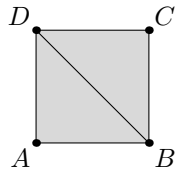
\epsfbox{fig3.1}$$
Here the segment $AB$ is identified with the segment
$DC$ and the segment $AD$ to the segment $BC$.

The integer areas of this two triangles are equal to $b+1$.

The integer lengths of the segment $BD$ is equal to $b+1$
and the integer of all others is equal to~1.

The integer angles $DAB$ and $DCB$ are equal to $b+1$,
and the others --- to~1.

All two triangles have the same integer affine type of the triangle
with the vertices $(0,0)$, $(0,1)$ and $b+1,0$ $1,0$
(b=5 at the picture):

$$
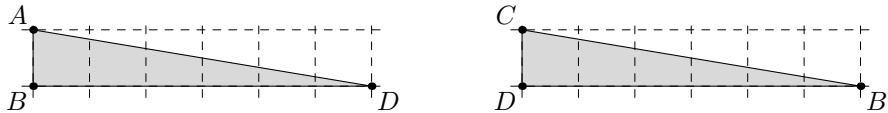
\epsfbox{fig3.6}$$

The integer distance from the triangle $ABD$ to the origin is equal to $a+2$,
and from the triangle $BCD$ is equal  to~1.
\end{example}

{\it Remark.} If we pose $b=0$, then we will have the sequence
of operators that have equivalent continued fractions with
the sequence of operators of Example~\ref{ex_k}.

\begin{statement}
The continued fractions are associated to the following operators
are equivalent (for integer $t\ge 0$):
$$
A_t=
\left(
\begin{array}{ccc}
0&1&0\\
0&0&1\\
1&1+t&-t-2\\
\end{array}
\right)
\quad
M_t=
\left(
\begin{array}{ccc}
0&0&1\\
1&0&-t-5\\
0&1&t+6\\
\end{array}
\right).
$$
\end{statement}

The proof of this statement is straightforward. It bases on the fact
that the operators $(I-A_t)^{-1}$ and $M_t$ are conjugate by
the operator $X$ in the group $SL(3,\z)$ (here $I$ --- is the identity
operator):
$$
X=
\left(
\begin{array}{ccc}
-1&-1&-2\\
0&0&-1\\
1&0&-1\\
\end{array}
\right).
$$
(i.e. $M_t=X^{-1}(I-A_t)^{-1}X$).

\begin{example}
The sequence of this example depends on an integer parameters $a \ge 1$.
The continued
fractions of the sequence are associated to the following operators:
$$
\left(
\begin{array}{ccc}
0&1&0\\
0&0&1\\
1&a&-2a-3\\
\end{array}
\right).
$$

The torus decomposition corresponding to this operator is homeomorphic
to the following one:
$$
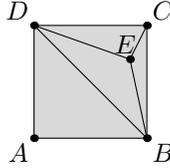
\epsfbox{fig3.7}$$
Here the segment $AB$ is identified with the segment
$DC$ and the segment $AD$ to the segment $BC$.

The integer areas, integer lengths and integer
angles of all four triangles are equal to $1$ here.

All two triangles have integer affine types of the simplest triangle:
$$
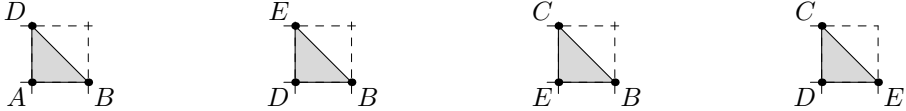
\epsfbox{fig3.8}$$

The integer distance from the triangle $ABD$ to the origin is equal to $a+2$,
from $BDE$~--- to $a+1$ and from $BCE$ and $CED$  --- to~1.
\end{example}

\begin{example}
The last sequence of this example that we will show in the issue
depends on two integer parameters $a >0$ and $b\ge 0$.
The continued
fractions of the sequence are associated to following operators:
$$
\left(
\begin{array}{ccc}
0&1&0\\
0&0&1\\
1&(a+2)(b+2)-3&3-(a+2)(b+3)\\
\end{array}
\right).
$$

The torus decomposition corresponding to this operator is homeomorphic
to the following one:
$$
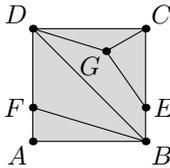
\epsfbox{fig3.9}$$
Here the segment $AB$ is identified with the segment
$DC$ and the polygonal line $AFD$ --- to the polygonal line $BEC$
(the point $F$ is identified with the point $E$).

The integer length of the segment $DB$ is equal to $b+2$
and all other integer lengths to~1.

The integer areas of the quadrangle $BECD$ is equal to $b+3$.
The integer areas of the triangle $BDF$ is equal to $b+2$,
and the integer lengths for all other triangles are equal to~1.

The integer angle $DFB$ is equal to $b+2$,
All other integer angles are equal to~1.

The triangles $ABF$, $CGE$ and $CDG$ have an integer affine type of the simplest
triangle;
the triangle $BDF$ has an integer affine type of the triangle
with the vertices $(0,0)$, $(b+2,0)$ and $(0,1)$;
the quadrangle $DBEC$ has an integer affine type of the quadrangle
with the vertices $(0,0)$, $(b+2,0)$, $(1,1)$ and $(0,1)$:

$$
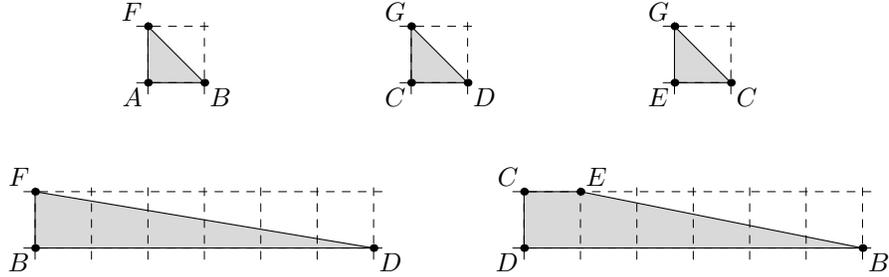
\epsfbox{fig3.10}$$

The integer distance from the triangles $ABF$ and $BFD$
to the origin is equal to~1,
from the triangle $CDG$ --- to $2+2a+2b+ab$,
from the triangle $CEG$ --- to $3+2a+2b+ab$,
from the quadrangle $DBEC$ --- to~1.
\end{example}

For more examples of such sequences and some method of
finding new such sequences, see~\cite{Kar1}.

\section{Some problems and conjectures.}

In the last section we formulate some of the most actual geometrical problems
and conjectures on Klein multidimensional continued fractions.
We start with the problem of the recognizing
the sails by their characteristics.

\begin{prob}\label{prob1}{\bf (Arnold)}
Find the collection of integer characteristics of the sail
that uniqu\-ely up to the equivalence relation determines the
two-dimensional $(n$-dimensional$)$ continued fraction.
\end{prob}

The same question is very important also for the periodic
two-dimensional sails of cubic irrationalities (and for dimension $n>2$).

For the one-dimensional case these characteristics are known.
Suppose we know the number of ordered sequences of integer lengths
and integer angles for one of the sails of the one-dimensional
continued fraction and their correspondence to each other.
This means that for any two segments with the common vertex
we know the integer lengths of this segments and the corresponding
integer angle between these edges.
The whole sail and even the whole continued fraction is
constructed by these characteristics (up to the equivalence relation)
in the unique way.
Moreover if we are given some arbitrary collection of this characteristics
then there exist the sail with this characteristics.
In the case of the periodic one-dimensional continued fraction it is
enough to know only the part of these two corresponding to each other
sequences of integer characteristics --- the integer lengths and angles
for some fundamental region (i.e. the decomposition of the circle).

There is no answer to this question in the two-dimensional case.
It is conjectured the following.

\begin{hyp}
To construct the whole continued fraction up to equivalence relation
it is sufficient to know only the following characteristics
for one of the sail:
the homeomorphic types of all facets for one of the sails
of the fraction, their adjacency to each other;
all their integer volumes and integer distances to them;
integer lengths of all edges;
and all integer angles for any two adjacent edges.
\end{hyp}

The whole classification of torus decompositions
than are obtained by some continued fraction of cubic irrationality
is unknown.
Here we formulate an important problem on the structure of torus
decomposition for the sails.

\begin{prob}{\bf (Arnold)}
Describe all torus decompositions that are possible for
some periodic two-dimensional continued fraction.
\end{prob}

Only a few is known in this direction.
There is a lot of trivial examples for torus decompositions
that do not corresponds to some sails. The following
example was constructed by the author in~\cite{Kar1}:

\begin{example}
The following torus decomposition does not correspond to any
sail for the periodic continued fraction of cubic irrationality
(there is one integer point in the interval $AD$, we denote it by $E$):
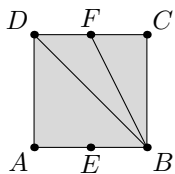
$$\epsfbox{fig2.1}$$
\end{example}

Here the polygonal line $AEB$ maps to the polygonal line $DFC$ under one of the
operators of the group $SL(3,\z)$ that preserves the sail.
As far $AEB$ is a segment, $DFC$ is also a segment. Therefore the points
$B$, $C$, $F$ and $D$ lie in the same plane. And hence $BF$ is not
a edge of some face.

Now we present the first nontrivial example of torus decomposition that was
made by E.~I.~Korkina in the work~\cite{Kor2}.

\begin{example}
Consider the simplest torus triangulation.
It consists of two triangles with the simplest
integer affine type of the triangle $(0,0)$, $(0,1)$ and $(1,0)$.
The integer distances to all of these faces is equal to~1.
This decomposition does not correspond to any
sail for the periodic continued fraction of cubic irrationality:
$$
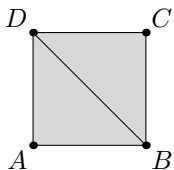
\epsfbox{fig3.1}$$
\end{example}

This example as far as a large number of sail calculations
allow us to formulate the following conjecture.

\begin{hyp}
For any torus decomposition for the continued fraction of the cubic
irrationality there exists some face with integer distance
to the origin more than one.
\end{hyp}

On the other hand each of the calculated sails have some face with
the integer distance to the origin equals one.

\begin{hyp}
For any torus decomposition for the continued fraction of the cubic
irrationality there exists some face with integer distance
to the origin equals one.
\end{hyp}

The following example is on the simplest torus decomposition,
that consist of one parallelogram.
The example was constructed by the author of this issue.

\begin{example}
Consider the torus decomposition consisting of the only one
face with integer affine type of the simplest parallelogram with the vertices
$(0,0)$, $(0,1)$, $(1,1)$ and $(1,0)$.
The integer distances to all of the faces is equal to some natural number
$n$.
This decomposition does not correspond to any
sail for the periodic continued fraction of cubic irrationality:
$$
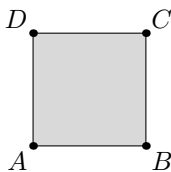
\epsfbox{fig2.2}$$
\end{example}

It seems that the torus decomposition with one rectangular face
is not possible for the decompositions associated to the
two-dimensional continued fractions of cubic irrationalities.
Moreover we conjecture the following.

\begin{hyp}
For any torus decomposition for the continued fraction of the cubic
irrationality there exists the face
that has an integer affine type of some triangle.
\end{hyp}

Let us describe the relation between two-dimensional
continued fractions of cubic irrationalities
and cubic extensions of the field of rational numbers.
For any such two-dimensional continued fraction
there exist an operator in $SL(3,\z)$ with irreducible characteristic
polynomial over the field $\q$ such that the continued fraction is associated
to the operator.
This polynomial determines the cubic extension of the field of
rational numbers.
All integer operators that identifies the continued fraction are all commutes
with $A$ and hence have the same cubic extension.
However the converse is not true (the example is from the work~\cite{Kar1}).

\begin{example}
The following two operators having the same characteristic
polynomial $x^3+11x^2-4x-1$ (and hence the same cubic extension)
define nonequivalent continued fractions:
$$
\left(
\begin{array}{ccc}
0&1&0\\
0&0&1\\
1&1&-2\\
\end{array}
\right)^3
,
\qquad
\left(
\begin{array}{ccc}
0&1&0\\
0&0&1\\
1&4&-11\\
\end{array}
\right).
$$
\end{example}

So, the natural question here is the following one.

\begin{prob}{\bf (Arnold)}
Classify the continued fractions that corresponds to the same cubic
extensions of the field of rational numbers.
\end{prob}

It is known almost nothing here. For example,
it is even unknown the finiteness of the number of possible
continued fractions associated to the same extension.
(For properties of cubic extensions of rational numbers
see in the work of B.~N.~Delone and D.~K.~Faddeev~\cite{Del1}.)

\begin{figure}
$$\epsfbox{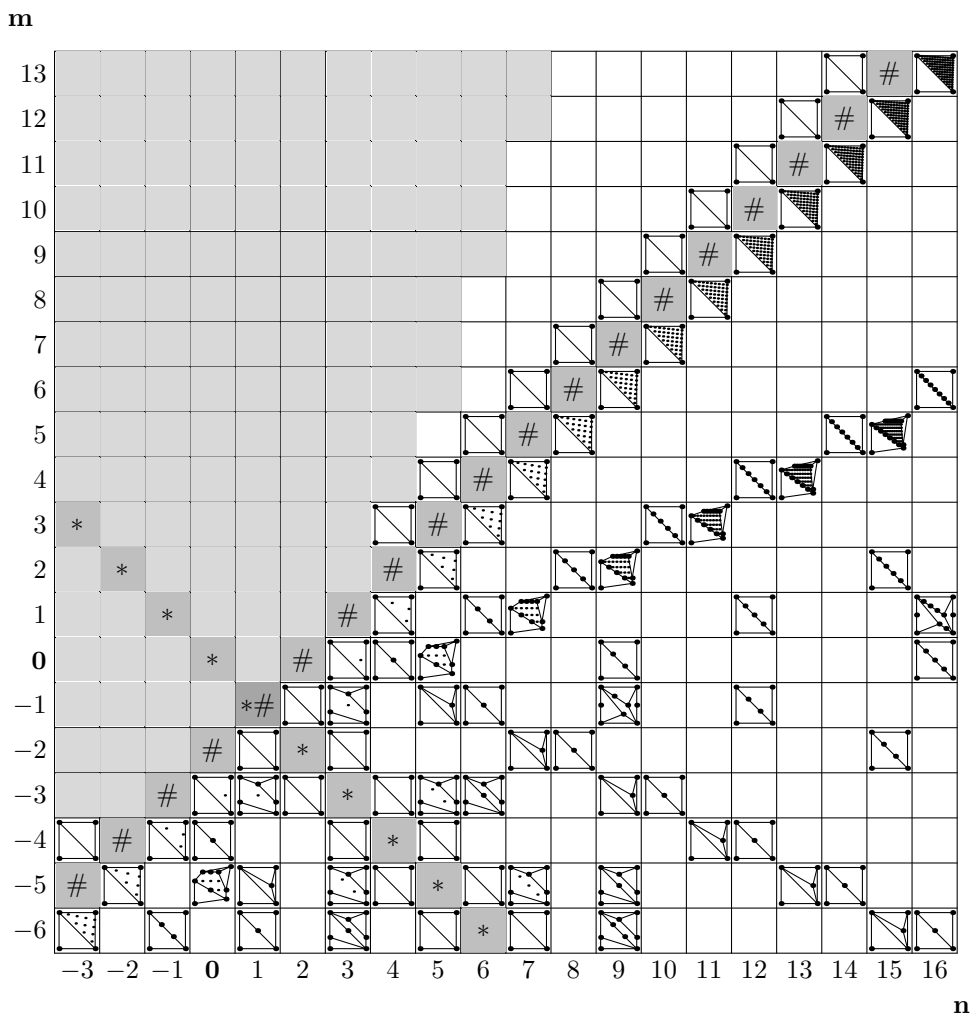}$$
\caption{Torus decompositions for operators $A_{m,n}$.}\label{pic2}
\end{figure}

In conclusion we show the table that was taken from my work~\cite{Kar2}
with squares filled with torus decomposition
of the sails for continued fractions associated to the family
of Sylvester operators
$$
A_{m,n}:=
\left(
\begin{array}{ccc}
0&1&0\\
0&0&1\\
1&-m&-n\\
\end{array}
\right)
,
$$
see fig~\ref{pic2}.
Some aspects of torus decomposition for the sail of the
two-dimensional continued fraction
for the cubic irrationality, associated to the operator $A_{m,n}$
is shown in the square sited at the intersection of the string with number $n$
and the column with number $m$.
If one of the roots of characteristic polynomial for the operator is equal to
1 or -1 at that than the square $(m,n)$ is marked with the sign $*$
or $\#$ correspondingly.
The squares that correspond to the operators which characteristic
polynomial has two complex conjugate roots
are painted over with light gray color.

\end{document}